\documentclass[a4paper]{article}
\usepackage{amssymb}
\usepackage{amsmath}
\usepackage{amsthm}
\usepackage{hyperref}
\usepackage{enumerate}
\usepackage{mathtools}
\usepackage{url}
\usepackage[T1]{fontenc}
\newcommand{\Cf}{{\mathbb C}}

\newcommand{\cA}{{\mathcal A}}

\newcommand{\cR}{{\mathcal R}}

\renewcommand{\phi}{\varphi} 

\title{Ratio bound (Lovász number) versus inertia bound}
\author{
 Ferdinand Ihringer
}
\date{13 May 2025}

\begin{document}
\maketitle

\begin{abstract}
  Matthew Kwan and Yuval Wigderson
  showed that for an infinite family of graphs, 
  the Lov\'asz number
  gives an upper bound
  of $O(n^{3/4})$ for the size of an independent set
  (where $n$ is the number of vertices),
  while the weighted inertia bound cannot do better than $\Omega(n)$.
  Here we point out that there is an infinite family of graphs for which 
  the Lovász number is $\Omega(n^{3/4})$,
  while the unweighted inertia bound is $O(n^{1/2})$.
\end{abstract}


\section{Introduction}

Initially, this note was written in 2023 as a small comment on a beautiful
result by Kwan and Wigderson \cite{KW2023}. Recently, Elphick, Tang, and Zhang \cite{ETZ2025}
did show that one can much better than the following example suggests.
This introduction got added for this occasion.

\section{An Example}
We mostly use the notation of \cite[\S2]{BCN} for association schemes.
Let $X$ be a finite set of size $n$. An {\it association scheme with $d$ classes}
is a pair $(X, \cR)$ such that
\begin{enumerate}[(i)]
 \item $\cR = \{ R_0, R_1, \ldots, R_d\}$ is a partition of $X \times X$,
 \item $R_0 = \{ (x, x): x \in X\}$,
 \item $R_i = R_i^T$, that is $(x, y) \in R_i$ implies $(y, x) \in R_i$,
 \item there are numbers $p_{ij}^k$ such that for any pair $(x, y) \in R_k$
 the number of $z$ with $(x, z) \in R_i$ and $(z, y) \in R_j$ equals $p_{ij}^k$.
\end{enumerate}
Note that some authors call $(X, \cR)$ as defined above a
{\it symmetric association scheme}.
For relations $R_i$, the $\{0,1\}$-adjacency matrices $A_i$ are defined by
\[
 (A_i)_{xy} = \begin{cases}
               1 & \text{ if } (x, y) \in R_i,\\
               0 & \text{otherwise.}
              \end{cases}
\]
As (i) holds, the matrices $A_i$ are linearly independent, and
as (iii) and (iv) hold, they generate a $(d+1)$-dimensional
commutative algebra $\cA$ of symmetric matrices, the {\it Bose-Mesner algebra}.
Since the $A_i$ commute, they can be diagonalized simultaneously
and we find a decomposition of $\Cf^n$ into a direct sum of
$d+1$ eigenspaces of dimension $f_j$ for $0 \leq j \leq d$.
As the all-ones matrix $J$ is in the span of $A_i$ and has
$n$ as an eigenvalue of multiplicity $1$, we may suppose that $f_0 = 1$.
If $\{ E_j: 0 \leq j \leq d \}$ is the basis of minimal idempotents of $\cA$,
then
\begin{align*}
 &f_j = \mathrm{rk}\, E_j = \mathrm{tr}\, E_j, && \sum_{j=0}^d E_j = I,
 && E_0 = n^{-1} J.
\end{align*}
Define matrices $P$ and $Q$ by
\begin{align*}
  &A_j = \sum_{i=0}^d P_{ij} E_i, && E_j = \frac{1}{n} \sum_{i=0}^d Q_{ij} A_i.
\end{align*}
Then $A_j E_i = P_{ij} E_i$ which shows that
the $P_{ij}$ are the eigenvalues of $A_j$.
Also note that $Q_{0j} = f_j$ as $\mathrm{tr}(E_j) = f_j$.

For a subset $Y$ of $X$ with characteristic vector $\chi$,
define a vector $a = (a_i)$, the {\it inner distribution of $Y$}, by
\[
 a_i := \frac{1}{|Y|} \chi^T A_i \chi = \frac{1}{|Y|} |\{ (x, y) \in Y \times Y \cap R_i \}|.
\]
Delsarte's linear programming bound states that
\[
 (aQ)_j \geq 0
\]
for all $0 \leq j \leq d$, see also Proposition 2.5.2 in \cite{BCN}.

We refer to \cite{KW2023} for details on the weighted ratio bound
and the weighted inertia bound (also called Cvetkovi\'c bound). It is well-known that
the weighted ratio bound (also called Hoffman bound) is a special case of the Lovász
number with equality in certain families of graphs.
For graphs which correspond to a union of relations in an
association scheme,
Delsarte's linear programming bound for independent sets
and the Lovász number are the same, see \cite{Schrijver1979}.
It is well-known that even the unweighted inertia bound
sometimes gives a better bound on the independence number
of a graph than the Lov\'asz number.
For instance, for the point graph of a generalized quadrangle
of order $(q, q^2)$, a graph with $(q^3+1)(q+1)$ vertices, 
the unweighted inertia bound is $q^3-q^2+q$, while the
Lov\'asz number is $q^3+1$.
Anurag Bishnoi asked if the inertia bound can also be 
asymptotically better than the Lov\'asz number
(as a parameter of the number of vertices $n$) \cite{AB2023}.
The purpose of this note is to point out that there 
exists a graph on $n$ vertices for which the 
Lov\'asz number is $\Omega(n^{3/4})$, but the
weighted inertia bound is $O(n^{1/2})$.
In \cite{CS1973} Cameron and Seidel describe a
$3$-class association scheme which has the following 
$P$- and $Q$-matrices (follow the instructions in
\cite[page 2]{DC2000} together with \cite{DCMM1997}
to obtain $P$ and $Q$ in a convenient manner):

\begin{align*}
  &P = \begin{pmatrix}1 & {{2}^{2 t}}{-}1 & {{2}^{4 t-2}}{+}{{2}^{3 t-2}}{-}{{2}^{2 t-1}}{-}{{2}^{t-1}} & {{2}^{4 t-2}}{-}{{2}^{3 t-2}}{-}{{2}^{2 t-1}}{+}{{2}^{t-1}}\\
1 & -1 & {{2}^{3 t-2}}-{{2}^{t-1}} & -{{2}^{3 t-2}}+{{2}^{t-1}}\\
1 & -1 & -{{2}^{t-1}} & {{2}^{t-1}}\\
1 & {{2}^{2 t}}-1 & -{{2}^{2 t-1}}-{{2}^{t-1}} & -{{2}^{2 t-1}}+{{2}^{t-1}}\end{pmatrix},\\
  &Q = \begin{pmatrix}1 & {{2}^{2 t}}-1 & {{2}^{4 t-1}}-3 \cdot {{2}^{2 t-1}}+1 & {{2}^{2 t-1}}-1\\
1 & -1 & -{{2}^{2 t-1}}+1 & {{2}^{2 t-1}}-1\\
1 & {{2}^{t}}-1 & -{{2}^{t}}+1 & -1\\
1 & -{{2}^{t}}-1 & {{2}^{t}}+1 & -1\end{pmatrix}
\end{align*}
Hence, using $f_j = Q_{0j}$, the graph with adjacency matrix $A_3$ has eigenvalues
\begin{itemize}
 \item ${{2}^{4 t-2}}-{{2}^{3 t-2}}-{{2}^{2 t-1}}+{{2}^{t-1}}$ with multiplicity $1$,
 \item $-{{2}^{3 t-2}}+{{2}^{t-1}}$ with multiplicity
 ${{2}^{2 t}}-1$,
 \item $2^{t-1}$ with multiplicity ${{2}^{4 t-1}}-3 \cdot {{2}^{2 t-1}}+1$,
 \item $-{{2}^{2 t-1}}+{{2}^{t-1}}$ with multiplicity ${{2}^{2 t-1}}-1$.
\end{itemize}
Hence, the unweighted inertia bound is
\[
  (2^{2t} -1) + (2^{2t-1} - 1) = 3 \cdot 2^{2t-1} - 2.
\]

The inner distribution $a$ of an independent set $Y$ of the graph has
the form $a = (1, x, y, 0)$, where $|Y| = 1+x+y$. Hence, the Lovász number is
the solution to the linear program which maximizes $1+x+y$
under the constraints that $(aQ)_j \geq 0$ for $0 \leq j \leq d$.
As $(aQ)_3 \geq 0$, we find that 
\[
  y \leq (x+1) \cdot (2^{2t-1}-1).
\]
As $(aQ)_2 \geq 0$, we find that 
\[
    (2^{2t-1}-1)x+(2^t-1)y \leq {{2}^{4 t-1}}-3 \cdot {{2}^{2 t-1}}+1.
\]
Clearly, $x = 2^{t}-1$ and $y=2^t \cdot (2^{2t-1}-1)$
maximizes $1+x+y$. As $(aQ)_j \geq 0$ for all $j$ for this solution,
this is an optimal solution. Hence, the
Lov\'asz number of the graph is $2^{3t-1}$.
Together with \cite{KW2023}, this shows the asymptotic incomparability of
the Lov\'asz number and the weighted inertia bound.

\noindent
Department of Mathematics,\\
Southern University of Science and Technology,\\ 
Shenzhen, China.\\
E-mail: {\tt ferdinand.ihringer@gmail.com}.

\end{document}